\documentclass[10pt,final]{siamltex}
\usepackage{graphicx}
\usepackage{latexsym}
\usepackage{amsmath,amssymb}
\usepackage{gensymb}
\usepackage{pstricks}
\begin{document}
\title{Mixed Maps for Kolmogoroff-Nagumo-Type Averaging on the Compact Stiefel Manifold}
\author{Simone Fiori, Tetsuya Kaneko, Toshihisa Tanaka %
\thanks{S. Fiori is with Dipartimento
        di Ingegneria dell'Informazione, Universit\`{a} Politecnica delle
        Marche, Via Brecce Bianche, I-60131 Ancona (Italy).\newline
        T. Kaneko and T. Tanaka are with the Department of Electrical and Electronic Engineering,
        Tokyo University of Agriculture and Technology (TUAT), 2-24-16, Nakacho,
        Koganei-shi,Tokyo 184-8588 (Japan).\newline
        This work was supported in part by KAKENHI 23300069.}}
\maketitle
\def\bbbr{{{\mathbb{R}}}}
\def\mdef{{\stackrel{{\mathrm{def}}}{=}}}
%
%
\begin{abstract}
The present research work proposes a new fast fixed-point averaging algorithm on the compact Stiefel manifold based on a mixed retraction/lifting pair. Numerical comparisons between fixed-point algorithms based on the proposed non-associated retraction/lifting map pair and two associated retraction/lifting pairs confirm that the averaging algorithm based on a combination of mixed maps  is remarkably less computationally demanding than the same averaging algorithm based on any of the constituent associated retraction/lifting pairs.
\end{abstract}
\normalsize
\section{Introduction}
The question about how to define the notion of \emph{mean value} over a set of structured samples has no intrinsic answer. Even in the simple case of a data set made of two positive real-valued numbers, there exist a variety of methods to define their mean value, each of which possesses different properties and leads to a different numerical result. The best known averages of two numbers $x,y>0$ are the arithmetic mean $\frac{1}{2}(x+y)$, the harmonic mean $\frac{2xy}{x+y}$ and the geometric mean $\sqrt{xy}$. A generalization of these is the Heinz mean $\frac{1}{2}(x^\alpha y^{1-\alpha}+y^\alpha x^{1-\alpha})$, with $0\leq \alpha\leq \frac{1}{2}$, that interpolates between the arithmetic ($\alpha=0$) and the geometric ($\alpha=\frac{1}{2}$) mean. Another averaging formula is given by the Heronian mean $\frac{1}{3}(x+\sqrt{xy}+y)$, that corresponds to a weighted sum of the arithmetic and the geometric means of the positive numbers $x$ and $y$. 
Most averaging formulas lead to instances of the \emph{Chisini mean} as defined in \cite{chipap}. The Chisini mean provides a good illustration of the fact that it is possible to define a mean value of a set of numbers without any particular requirement about, e.g., convexity. A function $f$ of $N$ real-valued variables leads to a Chisini mean value $\mu$ if, for every tuple $(x_1,x_2,\ldots,x_N)$, there exists a unique $\mu$ such that 
\begin{equation}\label{chimean}
 f(x_1,x_2,\ldots,x_N)=f(\mu,\mu,\ldots,\mu).
\end{equation}  
The definition of Chisini mean is really general and there were discussions whether it could be taken as a valid definition of mean (for a discussion, see, e.g., \cite{marichal}). But note that, for example, the arithmetic mean $\mu_\mathrm{a}\mdef\frac{1}{N}\sum_{k=1}^Nx_k$ of a tuple $(x_1,x_2,\ldots,x_N)$ is a Chisini mean with $f(x_1,x_2,\ldots,x_N)\mdef x_1+x_2+\cdots+x_N$. In fact, with this choice of the function $f$, the equation (\ref{chimean}) reads:
\begin{equation}
 \sum_{k=1}^Nx_k=\sum_{k=1}^N\mu=N\mu.
\end{equation}  
A rather general definition of mean value of a set of real-valued numbers led to the \emph{Kolmogoroff-Nagumo} (or \emph{quasi-arithmetic}) mean 
\begin{equation}\label{knm}
 \mu_\mathrm{KN}\mdef \varphi\left(\frac{1}{N}\sum_{k=1}^N \varphi^{-1}(x_k)\right)
\end{equation} 
for a continuous strictly monotonic function $\varphi$ (for a review, see \cite{marichal}). In the case of structured samples, such as structured matrices, the question about how to extend the definition of Kolmogoroff-Nagumo mean is rather involved. The Kolmogoroff-Nagumo mean turns out to be a special case of the Chisini mean. In fact, upon defining $f(x_1,x_2,\ldots,x_N)\mdef \varphi^{-1}(x_1)+\varphi^{-1}(x_2)+\cdots+\varphi^{-1}(x_N)$, it is readily verified that the mean $\mu_\mathrm{KN}$ is the unique solution to the Chisini equation (\ref{chimean})
\begin{equation}\label{knm1}
 \sum_{k=1}^N\varphi^{-1}(x_k)=\sum_{k=1}^N\varphi^{-1}(\mu)=N\varphi^{-1}(\mu).
\end{equation}  

Structured matrices appear in a variety of settings such as matrix-based optimization problems \cite{absilbook},  medical imaging \cite{bonnabel}, statistical analysis on manifolds \cite{Chikuse}, adaptive filtering \cite{fio:dsp}, computational ophthalmology \cite{siam11} and radio polarimetry \cite{hamaker}. The empirical mean is perhaps the most useful statistical characterization of a set of structured matrices \cite{arnaudon,bonnabel,fio:dsp,lim,turaga}, along with the empirical variance. When the matrices to average obey internal constraints such as orthogonality, the result of simple summation and division by the number of summands, which would represent an arithmetic mean, does not obey the same constraints, in general. Therefore, in order to obtain an empirical arithmetic mean of orthogonal matrices, it is necessary to establish a calculation method by considering the geometric structure of the matrix-space that the matrices to average belong to.  

To solve this problem, Kaneko, Tanaka and Fiori \cite{Kaneko} constructed an averaging algorithm on the compact Stiefel manifold which is a non-trivial extension of averaging algorithms on Lie-group-type manifolds presented in \cite{Fiori}. The underlying idea behind the algorithms developed in the contributions \cite{Kaneko} is that the sample-points on the Stiefel manifold are mapped onto a tangent space, where the average is taken, and then the obtained average point on the tangent space is brought back to the Stiefel manifold, via appropriate maps which are referred to as a {\it retraction map} and a {\it lifting map}. A variety of algorithms on the compact Stiefel manifold utilizing such maps were explained in  \cite{EAA}, whose convergence features and computational runtime were compared and which were applied to averaging real-world samples. 

The key point is that, in order to construct an averaging algorithm on the compact Stiefel manifold, it is necessary to construct an appropriate retraction/lifting pair. The averaging algorithms proposed in \cite{EAA} utilize QR-decomposition-based, Cayley-transform-based, polar-decomposition-based and orthographic-type retraction/lifting pairs. All these algorithms, however, suffer of a bloating of computational demand with the increase of the size of the processed matrices. For example, in the case of the polar-decomposition-based retraction/lifting pair, the retraction map may be expressed in closed form, while the computation of the associated lifting map requires solving a continuous-time algebraic Riccati equation.
 
In the present paper, in order to ease the computational demand of the previously-introduced class of averaging algorithms, a new Kolmogoroff-Nagumo-type averaging algorithm is proposed which exploits a combination of a closed form polar retraction map and a closed form orthographic lifting map. The implication of the proposed choice are evaluated analytically as well as numerically.
\section{Averaging algorithm based on a combination of retraction/lifting pairs}
The aim of the present section is to build-up an averaging algorithm on the compact Stiefel manifold based on the notion of \emph{mixed manifold retraction/lifting pair}. The algorithm, stemming from a non-linear matrix-type equation, is implemented by a fast fixed-point iteration scheme. 

The compact Stiefel manifold \cite{edelman} is defined by:
\begin{equation}
{\rm St}(p,n)\mdef\{X\in \mathbb{R}^{p\times n}|X^TX=I_n\},
\end{equation}
where symbol $I_n$ denotes a $n\times n$ identity matrix and $n\leq p$, namely, the manifold ${\rm St}(p,n)$ is the space of the `tall-skinny' orthogonal matrices. Its tangent space at a point $X\in{\rm St}(p,n)$ may be expressed as:
\begin{equation}
T_X{\rm St}(p,n)=\left\{V\in \mathbb{R}^{p\times n}|X^TV+V^TX=0\right\}.
\end{equation}
Each tangent space is a vector space under standard matrix addition and multiplication by a real scalar. 

The following measure of discrepancy $\delta$: St$(p, n)\times$ St$(p, n) \to \mathbb{R}_+^0$ between two Stiefel-manifold matrices is made use of:
\begin{equation}
\label{eq:delta}
\delta(X,Y)\overset{\rm def}{=}\| I_p-X^TY \|_{\rm F},\ X,Y\in{\rm St}(p,n),
\end{equation}
where symbol $\|\cdot\|_{\rm F}$ denotes the Frobenius norm.

{A \emph{retraction} at a point $X\in{\rm St}(p,n)$ of a Stiefel manifold is defined as a map $P_X:T_X{\rm St}(p,n)\rightarrow{\rm St}(p,n)$ (for a complete definition see, e.g., \cite{absilbook}), while a map $P_X^{-1}:{\rm St}(p,n)\rightarrow T_X{\rm St}(p,n)$ such that $P_X(P_X^{-1}(Q))=Q$, for $Q\in {\rm St}(p,n)$, is termed \emph{lifting} map (for a complete definition see, e.g., \cite{EAA}). Both maps are defined only locally, therefore, hereafter it is assumed that $X$ and $Q$ lay sufficiently close to each other to evaluate $P_X^{-1}(Q)$ and that $V$ is sufficiently close to $0$ to evaluate $P_X(V)$.}

Given a point $X\in{\rm St}(p,n)$ and a vector $V\in T_X{\rm St}(p,n)$, the \emph{polar-decomposition retraction} on the Stiefel manifold may be written in closed form \cite{absilbook}:
\begin{equation}\label{polret}
P_X(V) = (X+V)(I_n+V^TV)^{-\frac{1}{2}}.
\end{equation}
The associated \emph{polar-decomposition lifting} is denoted by $P_X^{-1}$.

The \emph{orthographic lifting map} \cite{malick} may be defined as follows:
\begin{equation}\label{eq1}
\hat{P}_X^{-1}(Q)=\pi_{T_X\mathrm{St}(p,n)}(Q-X),
\end{equation}
where $X,Q\in\mathrm{St}(p,n)$ and $\pi_{T_X\mathrm{St}(p,n)}:\bbbr^{p\times n}\rightarrow T_X\mathrm{St}(p,n)$ denotes a projection from the ambient space $\bbbr^{p\times n}$ into a tangent space $T_X\mathrm{St}(p,n)$. According to \cite{absilbook}, one such a projector is:
\begin{equation}\label{eq2}
\pi_{T_X\mathrm{St}(p,n)}(A)=(I_p-XX^T)A-X\mathrm{sk}(X^TA),
\end{equation}
where $A\in\bbbr^{p\times n}$ and $\mathrm{sk}(A)\mdef\frac{1}{2}(A^T-A)$. Plugging equation (\ref{eq2}) into equation (\ref{eq1}) gives:
\begin{eqnarray}
\nonumber
\hat{P}_X^{-1}(Q)&=&(I_p-XX^T)(Q-X)-X\mathrm{sk}(X^T(Q-X))\\
           &=&(I_p-XX^T)Q+\frac{1}{2}X(X^TQ-Q^TX).
\end{eqnarray}
The associated \emph{orthographic retraction map} is denoted by $\hat{P}_X$.
\subsection{Averaging method based on a mixed retraction/lifting pair}
The following steps lead to an equation characterizing the unknown empirical mean matrix $X\in{\rm St}(p,n)$, which represents an estimate of the actual center of mass $C\in{\rm St}(p,n)$ on the basis of the available information:
\begin{enumerate}
\item Map the points $X_k\in{\rm St}(p,n)$ belonging to a neighborhood of the sought-for mean-matrix $X\in{\rm St}(p,n)$ onto $T_X{\rm St}(p,n)$ by applying the lifting map $\hat{P}_X^{-1}$. Denote such points as $V_k\mdef \hat{P}^{-1}_X(X_k)$.
\item Compute the linear combination $\overline{V}=N^{-1}\sum^{N}_{k=1}V_k$. 
\item Bring back the mean vector $\overline{V}$ to ${\rm St}(p,n)$ by the retraction $P_X$ and get an empirical mean matrix $X=P_X(\overline{V})$.
\end{enumerate}
Summarizing the above procedure, a mean matrix $X\in{\rm St}(p,n)$
is the solution of the non-linear, matrix-type equation:
\begin{equation}\label{eq:mean}
X=P_X\left(\frac{1}{N}\sum^{N}_{k=1}\hat{P}^{-1}_X(X_k)\right)
\end{equation}
in the variable $X$. The analogy between the quasi-arithmetic mean for real numbers (\ref{knm}) and the equation (\ref{eq:mean}) is apparent, except that the equation (\ref{eq:mean}) is an implicit function of the mean (as both sides of the equation depend on the mean value) instead of being an explicit function as in the case of Kolmogoroff-Nagumo mean. 

The equation (\ref{eq:mean}) is solved numerically by
means of a fixed-point iteration algorithm, that generates a
sequence $X^{(i)}\in{\rm St}(p,n)$ of estimates of the
sought-for empirical mean matrix $X$, and that may be written as:
\begin{equation}\label{eq:mean_fp}
X^{(i+1)}=P_{X^{(i)}}\left(\frac{1}{N}\sum^{N}_{k=1}\hat{P}^{-1}_{X^{(i)}}(X_k)\right),\
i\geq 0,
\end{equation}
where matrix $X^{(0)}\in{\rm St}(p,n)$ denotes an initial guess. In the previous contributions \cite{Kaneko, EAA}, the algorithm (\ref{eq:mean_fp}) was based on associated retraction/lifting pairs, namely, either $P_X/P_X^{-1}$ or $\hat{P}_X/\hat{P}^{-1}_X$, while, in the present paper, the mixed retraction/lifting pair $P_X/\hat{P}^{-1}_X$ is made use of with the aim to decrease the computational burden of the algorithm.

Although, by definition of retraction/lifting pair, it holds $P_X\circ P^{-1}_X = \hat{P}_X\circ \hat{P}^{-1}_X=\mathrm{Id}_{\mathrm{St}(p,n)}$, it is recognized that the composition $P_X\circ \hat{P}^{-1}_X$ does not equal the identity map in $\mathrm{St}(p,n)$. In other terms, given matrices $X,Q\in \mathrm{St}(p,n)$, the discrepancy $\Delta_X(Q)\mdef\delta(P_X(\hat{P}^{-1}_X(Q)),Q)$ differs from zero, in general. Such a discrepancy may be evaluated as follows. First, note that:
\begin{equation}\label{discr}
 \Delta_X(Q)=\|I_n-Q^TP_X(\hat{P}^{-1}_X(Q))\|_{\mathrm{F}},
\end{equation}
define $M\mdef Q^TX$ and set $V=(I_p-XX^T)Q+\frac{1}{2}X(X^TQ-Q^TX)$. It holds that:
\begin{eqnarray}
V^TV=I_n+\frac{1}{4}(M^TM-M^{2T}-M^2-3MM^T),\\
Q^T(X+V)=I_n+M-\frac{1}{2}M(M+M^T),
\end{eqnarray} 
hence, the discrepancy (\ref{discr}) takes on the expression:
\begin{eqnarray}
\nonumber
\Delta&=&\left\|I_n-\left[I_n+M-\frac{1}{2}M(M+M^T)\right]\right.\times\\
      & &\left.\left[2I_n-\frac{1}{4}(M-M^T)^2-MM^T\right]^{-\frac{1}{2}}\right\|_{\mathrm{F}}.
\end{eqnarray}
Note that when $Q=X$ then $M=I_n$ and the above expression gives $\Delta=0$.
{ 
\subsection{Relationships with other contributions}
The iteration rule (\ref{eq:mean_fp}) generalizes the averaging methods proposed independently in \cite{EAA,R1}. The fixed-point solution $X$ of the iteration (\ref{eq:mean_fp}) satisfies the equation:
\begin{equation}\label{eq:zvf}
\sum^{N}_{k=1}\hat{P}^{-1}_X(X_k)=0,
\end{equation}
namely, the fixed point $X$ is a zero of the vector field $\mathcal{V}(X)\mdef\sum^{N}_{k=1}\hat{P}^{-1}_X(X_k)$. The Proposition 2.3 proved in \cite{R1} holds and ensures that, if the samples $X_k$ are not too spread (see \cite{R1} for a precise definition), the zero of the vector field $\mathcal{V}(X)$ is unique, hence, the iteration rule (\ref{eq:mean_fp}) is locally well-defined and leads to a unique average.

The fixed-point-type iteration rule (\ref{eq:mean_fp}) also generalizes the averaging methods based on the minimization of a spread function based on the Riemannian distance. A complete study of the convergence properties of such methods, based on a gradient-descent-type optimization procedure, is available in \cite{afsari}. An earlier study of an averaging method based on a Newton-type optimization procedure was proposed in \cite{hueper}.

The iteration rule (\ref{eq:mean_fp}) could be further generalized by introducing an unequal weighting scheme of the tangent-vectors $\hat{P}^{-1}_{X^{(i)}}(X_k)$ as suggested in the contribution \cite{R2} (see Algorithm 2). Such weighting scheme would lead to an iteration rule of the type:
\begin{equation}\label{eq:mean_wfp}
X^{(i+1)}=P_{X^{(i)}}\left(\frac{1}{N}\sum^{N}_{k=1}w_k^{(i)}\hat{P}^{-1}_{X^{(i)}}(X_k)\right),\
i\geq 0,
\end{equation}
where weights $w_k^{(i)}>0$ depend on the mutual distance from the samples and the current value of the iteration $X^{(i)}$.}
\section{Numerical experiments}
In the present section, the results of different experiments are illustrated to get an insight into the numerical behavior of the discussed retraction/lifting map pairs in the context of averaging over the compact Stiefel manifold. In the numerical experiments, the center $C \in {\rm St}(p,n)$ of the distribution of the samples $X_k$ is generated by computing the Q-factor of a thin-QR decomposition of a matrix randomly generated in $\mathbb{R}^{p\times n}$ with normally-distributed entries. The $N$ samples to average are generated by the rule $X_k=\exp(\sigma\Omega_k )C$, with $\Omega_k\overset{\rm def}{=}{\rm sk}(A_k)$, where $A_k$ is a matrix randomly generated in $\mathbb{R}^{p\times p}$ with normally-distributed entries, and $\sigma > 0$ controls the spread of the distribution \cite{Kaneko,EAA}. The initial guess $X^{(0)}$ in the fixed-point iteration algorithm (\ref{eq:mean_fp}) was chosen by slightly rotating the sample $X_1$ via a quasi-unit random rotation. {The numerical tests were performed by running a MATLAB\textsuperscript{\textregistered} 7 (64 bit) code on platform featuring an Intel\textsuperscript{\textregistered} Xeon\textsuperscript{\textregistered} (2.93 GHz) with 8 cores and 12GB RAM.}

\begin{figure}[t!]
\begin{center}
\includegraphics[width=3.8in]{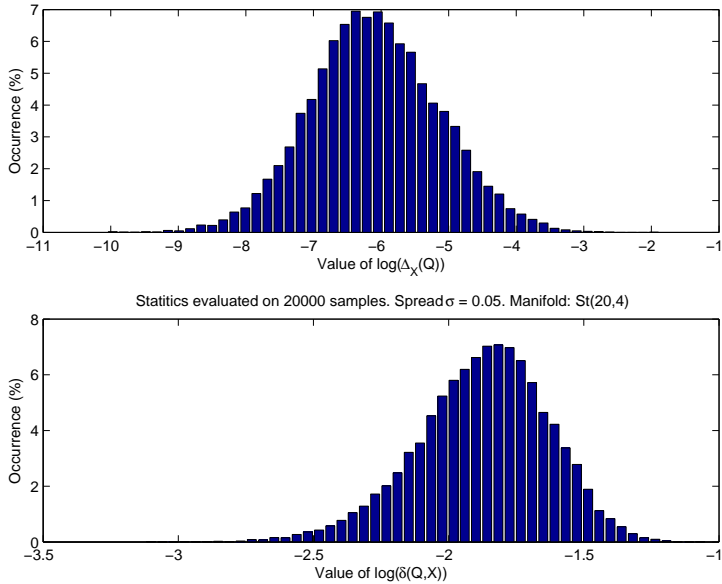}
\caption{Upper panel: Statistical distribution of the discrepancies $\Delta_X(Q)$ obtained by setting $X=C$ (center of the distribution) and $Q=X_k$ (random samples). Lower panel: Statistical distribution of the discrepancies $\delta(X,Q)$, again with $X=C$ and $Q=X_k$ (in logarithmic scales).}
\label{deltastat}
\end{center}
\end{figure}

\begin{figure}[t!]
\begin{center}
\includegraphics[width=3.8in]{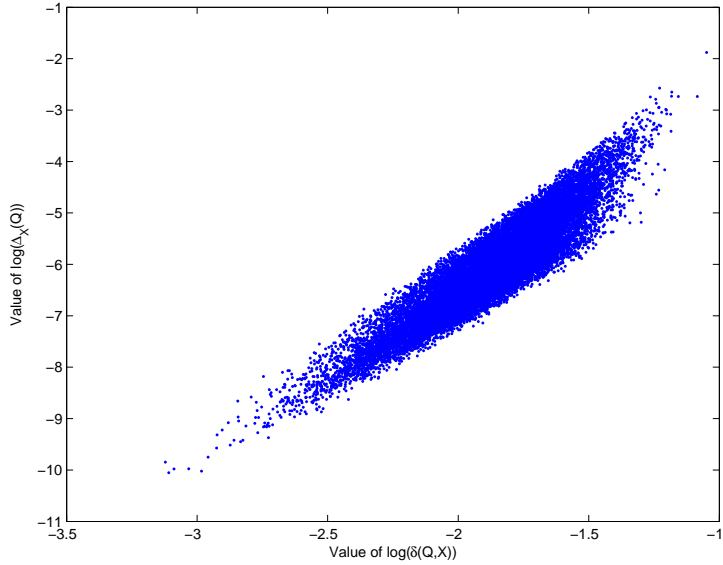}
\caption{Scatter plot of the values of the discrepancies $\Delta_X(Q)$ versus the values of the discrepancies $\delta(X,Q)$ obtained by setting $X=C$ and $Q=X_k$ (in logarithmic scales).}
\label{cloud}
\end{center}
\end{figure}

The first experiment aims at evaluating the numerical behavior of the compound map $P_X\circ \hat{P}^{-1}_X$. Figure~\ref{deltastat} illustrates the statistical distribution of the discrepancy $\Delta_C(X_k)$ compared to the distribution of the discrepancies $\delta(C,X_k)$. $N = 20,000$ samples $X_k\in{\rm St}(20,4)$ were generated by a spread factor of $\sigma = 0.05$. The discrepancy values $\Delta_C(X_k)$ distribute approximately around $10^{-6}$, confirming that the composition $P_X\circ \hat{P}^{-1}_X$, although not being an identity map, numerically behaves similarly to an identity map in the considered range. Figure~\ref{cloud} shows the scatter plot of the values of the discrepancies $\Delta_C(X_k)$ versus the values of the discrepancies $\delta(C,X_k)$. The relationship between them has a positive-correlation trend, showing that the more the arguments $Q$ and $X$ deviate from each other, the more the compound map $P_X\circ \hat{P}^{-1}_X$ deviates from the identity map.

\begin{figure}[t!]
\begin{center}
\includegraphics[width=3.8in]{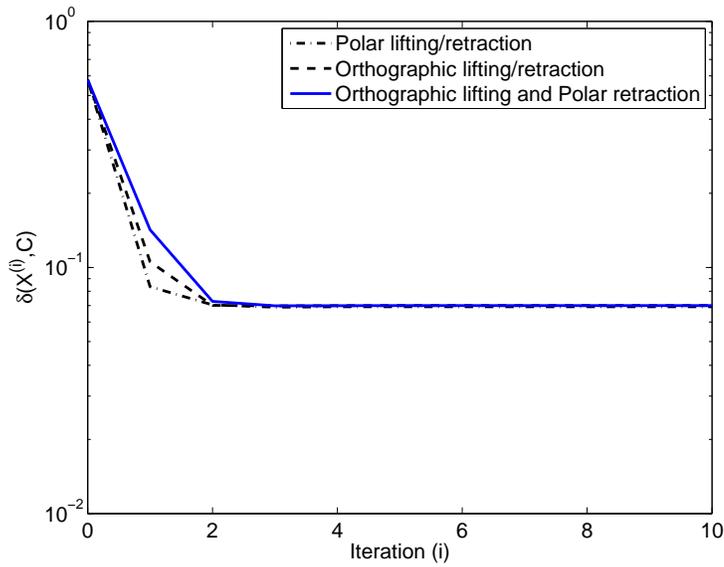}
\caption{Experiment about averaging on the manifold St(20,4). Index $\delta(X^{(i)}, C)$ during iteration. }
\label{conv1}
\end{center}
\end{figure}

The second experiment concerns the convergence properties of the fixed-point iteration algorithm (\ref{eq:mean_fp}) based on the polar retraction and orthographic lifting (namely, the mixed retraction/lifting pair $P_X/\hat{P}^{-1}_X$), the orthographic retraction/lifting pair ($\hat{P}_X/\hat{P}^{-1}_X$) and the polar retraction/lifting pair ($P_X/P^{-1}_X$). For this experiment, a number $N = 30$ of samples were generated on the manifold St$(20,4)$ with a spread parameter $\sigma = 0.2$. Figure~\ref{conv1} shows the values of the index $\delta(X^{(i)},C)$. The three algorithms behave satisfactorily and converge to solution-matrices that locate at similar distances to the actual center of the distribution. 

\begin{figure}[t!]
\begin{center}
\includegraphics[width=3.8in]{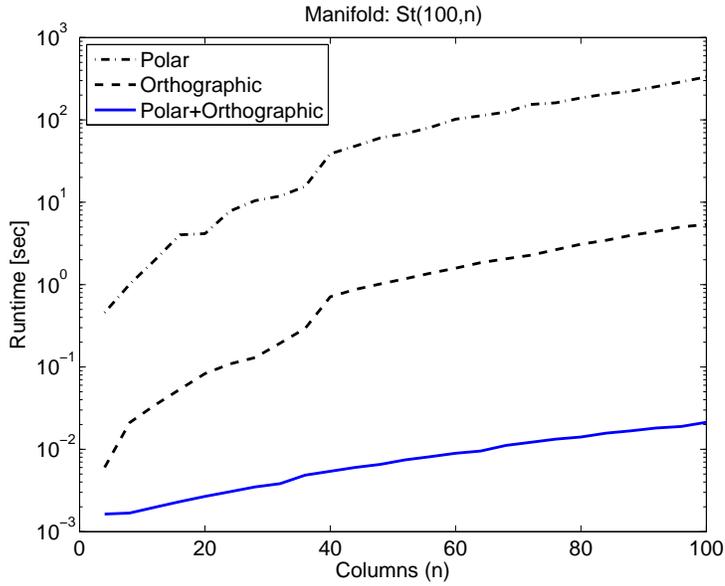}
\caption{Runtimes for the manifold St$(100,n)$ with the integer parameter $n$ varying. }
\label{runtime}
\end{center}
\end{figure}

\begin{figure}[t!]
\begin{center}
\includegraphics[width=3.8in]{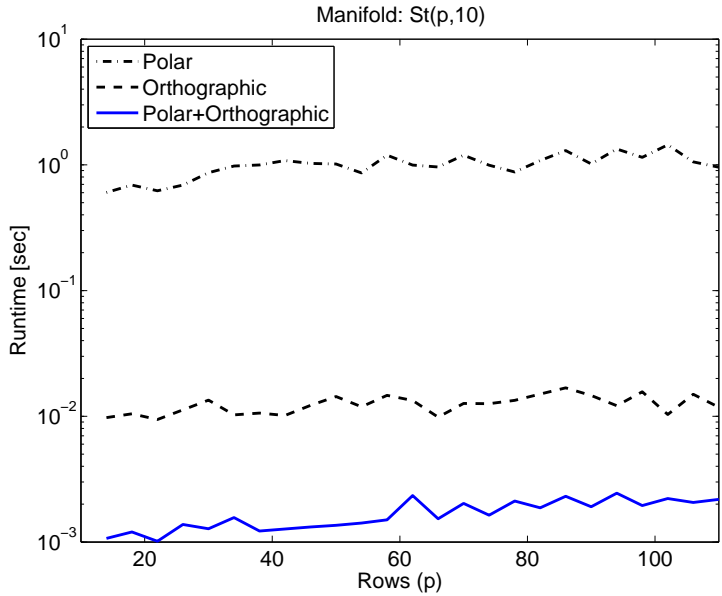}
\caption{Runtimes for the manifold St$(p,10)$ with the integer parameter $p$ varying. }
\label{runtime_rows}
\end{center}
\end{figure}

The third experiment aims at illustrating a comparison about the computational complexity of the three algorithms. A close inspection of the fixed-point algorithm (\ref{eq:mean_fp}) reveals that, for a general compact Stiefel manifold St$(p,n)$, the computational complexity is essentially a function of the number $n$. The Figure~\ref{runtime} shows the runtimes corresponding to the tested algorithms run on the manifold St$(100,n)$ with varying $n$. Such numerical simulation was performed with $N = 50$ samples generated with a spread-parameter value $\sigma = 0.01$. Each averaging experiment for each value of the index $n$ was repeated $100$ times to get rid of random fluctuations in the evaluation of runtimes. The obtained results indicate that the averaging algorithm based on a combination of the polar retraction map and the orthographic lifting map is much lighter than averaging methods based on associated retraction/lifting pairs in terms of computational burden. Likewise, Figure~\ref{runtime_rows} shows the runtimes corresponding to the tested algorithms run on the manifold St$(p,10)$ with varying $p$. Such numerical simulation was performed with $N = 50$ samples generated with a spread-parameter value $\sigma = 0.01$. The experiment for each $p$ was repeated on $100$ independent trials to get rid of random fluctuations. The obtained results indicate the little dependence of the computational complexity of the algorithm from the dimension $p$.
\section{Conclusion}
The present research work proposed a new averaging algorithm on the compact Stiefel manifold based on a combination of non-associated (mixed) retraction and lifting maps. The combined retraction/lifting pair was studied both analytically and numerically and the new averaging algorithm was tested numerically. The obtained numerical results confirm that the composition $P_X\circ \hat{P}^{-1}_X$, although not being an identity map, numerically behaves like an identity map in the considered range and that the new algorithm behaves satisfactorily and converges to solution-matrices that deviate from the actual center of the distribution of a similar amount. A comparison between three fixed-point algorithms based on the non-associated polar retraction and orthographic lifting, on the associated orthographic retraction/lifting pair and on the associated polar retraction/lifting pair reveals that the averaging algorithm based on a combination of the polar retraction map and the orthographic lifting map is remarkably superior to the averaging algorithms based on associated retraction/lifting pairs in terms of computational demand, resulting the lightest one.

%

%

\begin{thebibliography}{99}
%
\bibitem{absilbook}
P.-A. Absil, R. Mahony and R. Sepulchre, {\it Optimization
Algorithms on Matrix Manifolds.} New Jersey: Princeton University
Press, 2008.
%
\bibitem{malick}
P.-A. Absil and J. Malick, ``Projection-like retractions on matrix
manifolds,'' {\it SIAM Journal on Optimization}, Vol. 22. No. 1, pp.
135 -- 158, 2012.
%
\bibitem{afsari}
B. Afsari, R. Tron and R. Vidal, ``On the convergence of gradient descent for finding the Riemannian center of mass,'' {\it SIAM Journal on Control and Optimization}. Submitted
%
\bibitem{arnaudon}
M. Arnaudon, F. Barbaresco and L. Yang, ``Medians and means in
Riemannian geometry: Existence, uniqueness and computation,'' In
\emph{Matrix Information Geometries} (R. Bhatia and F. Nielsen,
Eds.), Springer, 2012.
%
\bibitem{bonnabel}
S. Bonnabel and R. Sepulchre, ``Riemannian metric and geometric mean
for positive semidefinite matrices of fixed rank,''  {\it SIAM Journal 
on Matrix Analysis and Applications}, Vol. 31, No. 3, pp. 1055 -- 1070,
August 2009. 
%
\bibitem{Chikuse} 
Y. Chikuse, \emph{Statistics on Special Manifolds}. Springer, 2003.
%
\bibitem{chipap}
O. Chisini, ``Sul concetto di media,'' {\it Periodico di Matematiche}, 
Vol. 4, pp. 106 -- 116, 1929.
%
\bibitem{edelman}
A. Edelman, T. A. Arias, and S. T. Smith, ``The geometry of
algorithms with orthogonality constraints," {\it SIAM Journal on
Matrix Analysis and Applications}, Vol. 20, No. 2, pp. 303 -- 353,
1998.
%
\bibitem{fio:dsp}
S. Fiori, ``On vector averaging over the unit hyphersphere," {\it Digital
Signal Processing}, Vol. 19, No. 4, pp. 715 -- 725, July
2009.
%
%
\bibitem{siam11} S. Fiori, ``Solving Minimal-Distance Problems over the Manifold
of Real Symplectic Matrices,'' {\it SIAM Journal on Matrix Analysis
and Applications}, Vol. 32, No. 3, pp. 938 -- 968, 2011.
%
\bibitem{Fiori} S. Fiori and T. Tanaka, ``An algorithm to compute averages on matrix Lie groups," {\it IEEE Transactions on Signal Processing}, Vol. 57, No. 12, pp. 4734 -- 4743, December 2009. 
%
\bibitem{R1} P. Grohs, ``Geometric multiscale decompositions of dynamic low-rank matrices,'' SAM Report 2012-03, ETH Zurich, February 2012 
%
\bibitem{hamaker}
J.P. Hamaker, ``Understanding Radio Polarimetry IV: The full-coherency 
analogue of scalar self calibration,'' {\it Astronomy and Astrophysics 
Supplement Series}, Vol. 143, No. 3, pp. 515 -- 534, May 2000.
%
\bibitem{hueper}
K. H\"{u}per and J. Trumpf, ``Newton-like methods for numerical optimization on manifolds,'' in {\it Conference Record of the Thirty-Eighth Asilomar Conference on Signals, Systems and Computers} (Asilomar Hotel and Conference Grounds, November 7-10, 2004), Vol. 1, pp. 136 -- 139, 2004.
%
%
\bibitem{Kaneko} T. Kaneko, T. Tanaka, and S. Fiori, ``A method to compute averages over the compact Stiefel manifold," in {\it Proceedings of the IEEE International Conference on Acoustic, Speech and Signal Processing} (ICASSP 2012, Kyoto, Japan, March 25 - 30, 2012), pp. 3829 -- 3832, 2012.
%
\bibitem{EAA} T. Kaneko, S. Fiori and T. Tanaka, ``Empirical arithmetic averaging over the compact Stiefel manifold," 
{\it IEEE Transactions on Signal Processing}, Vol. 61, No. 4, pp. 883 -- 894, February 2013.
%
%
\bibitem{lim}
J.D. Lawson and Y. Lim, ``The geometric mean, matrices, metrics, and
more,'' {\it The American Mathematical Monthly}, Vol. 108, No. 9,
pp. 797 -- 812, November 2001.
%
\bibitem{R2} {J. Marks, M. Kirby and C. Peterson, ``A normal/tangent bundle algorithm for representing point clouds on Grassmann and Stiefel manifolds," Available online at \verb+http://www.mendeley.com/download/public/11078691/+ \verb+4745581302/48b1f87dfc3dd21a7140d0ae524963a79f2721+ \verb+d3/dl.pdf+}
%
\bibitem{marichal}
J.-L. Marichal, ``On an axiomatization of the quasi-arithmetic mean values without the symmetry axiom,'' {\it Aequationes Mathematicae}, Vol. 59, No. 1-2, pp. 74 -- 83, 2000.
%
%
\bibitem{turaga} P. Turaga, A. Veeraraghavan, A. Srivastava and R.
Chellappa, ``Statistical computations on Grassmann and Stiefel manifolds 
for image and video-based recognition,'' {\it IEEE Transactions on Pattern
Analysis and Machine Intelligence}, Vol. 33, No. 11, pp. 2273 -- 2286, 2011.
%
\end{thebibliography}
\end{document}